\documentclass[times,sort&compress,3p,12pt]{elsarticle}
\usepackage{natbib}
\usepackage{setspace}
\usepackage[labelfont=bf]{caption}

\usepackage{amssymb}
\usepackage[utf8]{inputenc}
\usepackage{amssymb,amsmath,mathrsfs}
\usepackage{amsthm}
\usepackage[colorlinks,
            linkcolor=blue,
            anchorcolor=red,
            citecolor=green
            ]{hyperref}
\numberwithin{equation}{section}
\usepackage{graphicx}
\usepackage{amsmath}
\usepackage{amsmath,amsfonts,amssymb,amsthm,booktabs,color,epsfig,graphicx,hyperref,url}
\usepackage{enumitem}
\usepackage{array}
\usepackage{graphicx}
\usepackage{tikz}
\usetikzlibrary{backgrounds}
\usepackage{xcolor}
\usepackage{yhmath} % 引入自适应宽度的 hat 支持
\usepackage{dsfont} % 在导言区引入宏包

\theoremstyle{plain}% Theorem-like structures provided by amsthm.sty

\definecolor{mycustompurple}{RGB}{154, 36, 79}
\usepackage[utf8]{inputenc}

\usepackage{hyperref} %

\hypersetup{
	colorlinks=true,            %
	linkcolor=blue,              %
	citecolor=blue,            %
	filecolor=mycustompurple,   %
	urlcolor=magenta            %
}

\newtheorem{theorem}{Theorem}[section]
\newtheorem{proposition}[theorem]{Proposition}

\newtheorem{lemma}[theorem]{Lemma}

\newtheorem{remark}[theorem]{Remark}

\def\C{\mathbb C}
\def\R{\mathbb R}
\def\Re{{\rm Re\hspace{0.05cm}}}
\def\Im{{\rm Im\hspace{0.05cm}}}
\def\d{\hspace{0.05cm}\mathrm{d}}

\def\Y{Y_{r,u}}

\begin{document}
	
	\begin{frontmatter}
		
		\title	{
			Persistence in weighted Sobolev spaces
for the Calogero–Moser derivative nonlinear Schr\"odinger equation}
		
		\author[1]{Ruoci Sun}
		\author[2]{Yi Zhao\texorpdfstring{\corref{mycorrespondingauthor}}{}}

		\address[1]{Laboratory of Mathematics and Complex Systems (Ministry of Education), 
			School of Mathematical Sciences, 
			Beijing Normal University, 
			Beijing 100875, China.}
		
		\address[2]{College of Science, Minzu University of China, Beijing 100081, China.}
		
		\cortext[mycorrespondingauthor]{Corresponding author: \url{zhaoyi202406@163.com}}  % 补充邮箱
		
		\begin{abstract}

		This paper is devoted to the study of  invariant weighted Sobolev spaces for the Calogero--Moser
derivative nonlinear Schr\"odinger (CMDNLS) equation on the real line.
For \(s\geq 1\) and \(0\leq r\leq s\), we first prove that the
\(H^s(\mathbb R)\) flow preserves the condition
\((1+x^2)^{r/2}u\in L^2(\mathbb R)\).
The proof uses truncated weighted energy estimates and weighted
interpolation. 
We then establish the persistence of weighted regularity for spatial
derivatives. For integers \(p\geq 3\) and \(1\leq m\leq p-2\), and
for \(0\leq r\leq p-m \), the \(H^p(\mathbb R)\) flow
preserves the conditions
\((1+x^2)^{r/2}\partial_x^j u\in L^2(\mathbb R)\) for
\(0\leq j\leq m\).
This result combines weighted energy estimates, the unitary equivalence and  the integrable
structure of the CMDNLS equation.

		\end{abstract}

		
	\end{frontmatter}
\tableofcontents

\section{Introduction}

The \emph{Calogero--Moser derivative nonlinear Schr\"odinger equation} (CMDNLS), 
also known as the \emph{continuum Calogero--Moser equation} (CCM), is given by
\begin{equation}\label{foc_de_CMDNLS}
  \partial_t u = i\,\partial_x^2 u + \sigma i\, u \,\mathbb{H}\partial_x(|u|^2) 
  + \sigma\, u \,\partial_x(|u|^2), \qquad (t,x)\in \mathbb{R}\times\mathbb{R},
\end{equation}
where $u=u(t,x)$ is a complex-valued function, $\sigma=\pm 1$, 
and $\mathbb{H}$ denotes the Hilbert transform
\begin{equation}\label{eq:Hilbert}
  \mathbb{H}f(x) = \frac{1}{\pi}\,\mathrm{p.v.}\int_{\mathbb{R}} \frac{f(y)}{x-y}\,\mathrm{d}y.
\end{equation}
Here $\mathrm{p.v.}$ stands for the Cauchy principal value. 
The case $\sigma=1$ is referred to as the \emph{focusing} CMDNLS equation
\begin{equation} \tag{$\rm{CMDNLS}^+$} \label{foc_CMDNLS}
  \partial_t u = i\,\partial_x^2 u +  i\, u \,\mathbb{H}\partial_x(|u|^2) 
  + \, u \,\partial_x(|u|^2), \qquad (t,x)\in \mathbb{R}\times\mathbb{R},
\end{equation}
and  $\sigma=-1$ corresponds to the \emph{defocusing} one

\begin{equation} \tag{$\rm{CMDNLS}^-$} \label{de_CMDNLS}
  \partial_t u = i\,\partial_x^2 u -i\, u \,\mathbb{H}\partial_x(|u|^2) 
  - \, u \,\partial_x(|u|^2), \qquad (t,x)\in \mathbb{R}\times\mathbb{R}.
\end{equation}

The focusing CMDNLS equation \eqref{foc_CMDNLS} arises as the thermodynamic 
continuum limit of the classical Calogero--Moser Hamiltonian system \cite{Abanov}
\begin{equation}\label{eq:CM}
  \mathcal{H}_{\mathrm{CM}} = \frac{1}{2}\sum_{j=1}^{N} p_j^2 
  + \sum_{1\le j\neq k\le N}\frac{1}{2(x_j-x_k)^2},
\end{equation}
which describes a completely integrable system of $N$ particles on the line 
with pairwise inverse-square interactions. 
For sufficiently large positive times, pure $N$-soliton solutions
of the focusing CMDNLS equation admit rational representations
with $N$ distinct poles $z_1(t),\ldots,z_N(t)$ in the lower
half-plane. These poles satisfy the complexified Calogero--Moser
system
\cite{GerardCPAM}
\begin{equation}\label{eq:complex-CM}
  \partial_t^2{z}_j(t) = -\sum_{k\neq j}\frac{8}{(z_k-z_j)^3}, 
  \qquad j=1,\dots,N,
\end{equation}

  The defocusing CMDNLS equation \eqref{de_CMDNLS} appears as the deep-water 
limit ($\delta\to\infty$) of the \emph{intermediate nonlinear Schr\"odinger 
equation} (INLS) \cite{Peli-PLA,Peli-JMP}
\begin{equation}\label{eq:INLS}
  i\,\partial_t u + \partial_x^2 u + u\,(i+T_\delta)\partial_x(|u|^2)=0,
\end{equation}
where the singular integral operator $T_\delta$ is defined by
\begin{equation}\label{eq:Tdelta}
  T_\delta f(x) := \mathrm{p.v.}\int_{\mathbb{R}} 
  \coth\!\left(\frac{\pi(y-x)}{2\delta}\right) f(y)\,\frac{\mathrm{d}y}{2\delta},
  \qquad \delta>0.
\end{equation}
The INLS equation \eqref{eq:INLS} interpolates between two fundamental models. 
After a suitable rescaling, its formal shallow-water limit
($\delta\to0$) yields the classical nonlinear Schr\"odinger
equation, whereas its deep-water limit gives
\eqref{de_CMDNLS}.  This mirrors the well-known relationship between the Korteweg--de Vries (KdV) 
equation, the intermediate long wave (ILW) equation, and the Benjamin--Ono (BO) 
equation --- where KdV and BO arise respectively as the shallow-water and 
deep-water limits of ILW.

Throughout this paper, the Fourier transform is denoted by $\hat f$ with
$$
\hat f(\xi):=\int_{\mathbb{R}} f(x) e^{-i\xi x} \d x, \quad\forall f\in L^1=L^1(\mathbb{R}; \mathbb{C}).
$$  
For  $s>0$, the Sobolev space  $H^s=H^s(\R; \C)$ is given by  
 
\[
H^s:=\left\{
f\in L^2=L^2(\R; \C): \|f\|_{H^s}^2 := \frac{1}{ 2 \pi} \int_{\mathbb{R}} (1+\xi^2)^s  |\hat f(\xi)|^2 \mathrm{d} \xi   <+\infty
\right\},
\]
 
The focusing CMDNLS equation is globally well-posed in the Hardy-Sobolev space 
\begin{equation}\label{HardySobolevHs+}
     H^s_+:= \{f \in H^s: \mathrm{supp}\hat{f} \subset [0, +\infty) \}, \quad s\in \R
 \end{equation}
for $s\ge 1$  with  initial data  $\| u_0\|^2_{L^2}\le  2\pi$ \cite{GerardCPAM}. 
Let $L^2_+=L^2_+(\R; \C):=H^0_+$ denote the space that consists of all $L^2$ functions whose  Fourier transforms are supported in $[0,\infty)$. The Cauchy-Szeg\H{o} projection is denoted by $\Pi=\frac{\rm{id}+i\mathbb{H}}{2}: L^2 \to L^2_+$. So for any $f\in L^2$, we have 
\[
\widehat{\Pi f}(\xi)=\mathds{1}_{\xi\geq 0}\hat f(\xi)=\begin{cases} \hat{f}(\xi), & \xi \geq 0, \\ 0, & \xi < 0. \end{cases}
\]

The focusing equation is $L^2$-critical and admits the phase,
scaling, and spatial translation symmetries
\[
u(t,x)\mapsto e^{i\theta}\lambda^{1/2}
u(\lambda^2t,\lambda x+x_0),\qquad
\theta,x_0\in\mathbb R,\quad\lambda>0.
\]
The equation on the full space also admits the Galilean transformation
\[
u(t,x)\mapsto e^{i\eta x-it\eta^2}u(t,x-2t\eta),
\qquad\eta\in\mathbb R.
\]
This transformation preserves the nonnegative Fourier support of
every Hardy function when $\eta\ge0$.
The focusing CMDNLS equation admits infinitely many conservation laws. The three most fundamental ones are 
\begin{equation}
    E(u)=\frac 12 \int_\R |\partial_x u - i \Pi(|u|^2)u|^2 \d x,\quad M(u)=\int_\R |u|^2 \d x, \quad P(u)=\Re \int_\R \big( 
    -i\overline u \partial_x u-\frac 12 |u|^4 \big) \d x.
\end{equation}

Eq. \eqref{foc_de_CMDNLS} is completely integrable. The Lax pair structure is introduced in \cite{GerardCPAM}.  
For $u\in H^1_+$, the Lax operator $(\mathbb{L}_u, \rm{Dom}(\mathbb{L}_u)=H^1_+)$ is defined by 
\begin{equation}
    \mathbb{L}_u(h)=-i\partial_x h-\sigma u\Pi(\overline u h ), \quad \forall h\in H^1_+.
\end{equation}
The operator $(\mathbb{L}_u, \rm{Dom}(\mathbb{L}_u)=H^1_+)$ is an unbounded self-adjoint operator on the Hilbert space $L^2_+$.
If $u$ solves \eqref{foc_de_CMDNLS}, then there exists an operator $\mathbb{B}_u$ depending on $u$ such that 
\begin{equation}
    \partial_t \mathbb{L}_u=[\mathbb{B}_u, \mathbb{L}_u],
\end{equation}
where $[\cdot,\cdot]$ denotes the commutator.
Using the Lax pair structure, the direct scattering theory  of \eqref{foc_de_CMDNLS} is studied in \cite{Frank}\cite{arxiv}.
Weighted Sobolev regularity is relevant to the time evolution of the Jost functions in the direct scattering
problem for the CMDNLS equation.
For  any integer $s\ge 0$ and  $r\ge0$, define
\begin{equation}\label{H_r^s}
H_r^s=\left\{f\in H^s\,\middle|\,
\langle x\rangle^r\partial_x^j f\in L^2,
\quad 0\le j\le s\right\},\qquad
\|f\|_{H_r^s}^2=
\sum_{j=0}^s\|\langle x\rangle^r\partial_x^j f\|_{L^2}^2,
\end{equation}
where $\langle x\rangle=(1+x^2)^{1/2}$. In particular, $H_r^0=L_r^2$.
We study whether this weighted regularity persists throughout
the lifespan of a sufficiently regular solution.
In particular, we show that
$H^4\cap H_1^2$ is invariant under the $H^4$ flow
of \eqref{foc_de_CMDNLS} throughout its lifespan \cite{arxiv}. Weighted derivative persistence was established in
\cite[Theorem 1.7]{arxiv} for $H^4_+\cap H_1^2$
under the focusing CMDNLS flow.
Theorem~\ref{H^N_1} removes the Hardy space restriction,
treats both signs of the equation, and allows general
derivative and weight. This is also named the persistence space, i.e., the solution $u(t)$ belongs to the same space as the initial data $u_0$ for any  $t$ in the lifespan.

Persistence means that a solution retains the regularity and spatial
decay of its initial datum throughout its lifespan. The persistence property is an essential ingredient of well-posedness: together with existence, uniqueness and continuous dependence. It guarantees that the flow of the equation leaves the chosen function space invariant, so that the solution map can be regarded as a self-map on that space. It shows that the smoothness of the initial data is preserved within the lifespan, which serves as a basis for establishing blow-up  and global existence results. Moreover, persistence in weighted spaces yields information on the spatial decay of solutions and constitutes a key tool in the study of unique continuation properties and asymptotic behavior. It also underlies the long-time dynamics, such as stability analysis and scattering theory.

\subsection{Main results}

The principal goal of the paper is to establish the persistence of regularity and decay for solutions of \eqref{foc_de_CMDNLS}. More precisely, we seek function spaces that are invariant under the flow of \eqref{foc_de_CMDNLS}.  To this end, we pose the initial data in the weighted $L^2$ space 
\[
L_r^2=\left\{
f:\R\to \C: \langle x \rangle^r f\in L^2
\right\},\quad \langle x \rangle : =\sqrt{1+x^2}, \quad r\ge 0. 
\]

Our first result concerns  the invariance of the weighted Sobolev space $H^s\cap L_r^2$ for $s\ge 1$ and $0\le r \le s$ under the CMDNLS flow. 
First, recall the local well-posedness result  of \eqref{foc_de_CMDNLS} on the full Sobolev space \cite{laurens2025}. 
 Let $u_0\in H^s$ and $s>\frac 14$. Then there exists a time 
\begin{equation}
\mathfrak{T}=\mathfrak{T}(\| u_0\|_{H^s})>0,
\end{equation}
and a unique solution $u\in C([-\mathfrak{T}, \mathfrak{T}]; H^s)$ of \eqref{foc_de_CMDNLS}  with the initial datum $u_0$.
\begin{theorem}\label{foc_invariant_HsLr}
Let $s\ge1$, $0\le r\le s$, and $0<T\le\mathfrak T$.
Suppose that $u\in C([-T,T];H^s)$ solves \eqref{foc_de_CMDNLS}
with $u(0)=u_0\in H^s\cap L_r^2$. Then
\[
u\in C([-T,T];H^s\cap L_r^2).
\]
\end{theorem}

\begin{remark}
The persistence property here differs from the Benjamin--Ono equation's
\begin{equation}
    \partial_t v+v\partial_{x}v+\mathbb{H}\partial_{x}^2v=0, \quad (t,x)\in \R\times \R.
\end{equation}
    The zero-frequency jump used in  Benjamin--Ono argument (see the proof of Theorem 2 and Theorem 3 in \cite{Ponce2011})  is absent in \eqref{foc_de_CMDNLS}. Because there exists a term
    \[
\partial_\xi^2 e^{-it\xi|\xi|}
=\bigl(-4t^2\xi^2-2it\,\operatorname{sgn}(\xi)\bigr)
  e^{-it\xi|\xi|},
\]
  in  the linear part of Benjamin--Ono equation. However, the linear part of \eqref{foc_de_CMDNLS} involves 
    \[
\widehat{e^{it\partial_x^2}f}(\xi)
=e^{-it\xi^2}\widehat f(\xi),
\qquad
\partial_\xi^2 e^{-it\xi^2}
=\bigl(-4t^2\xi^2-2it\bigr)e^{-it\xi^2}.
\] 
The sign discontinuity involved in the Benjamin--Ono analysis
\cite{Ponce2011} is therefore absent from the linear Schr\"odinger
phase in equation \eqref{foc_de_CMDNLS}. 
   The following construction gives nonzero solutions of the CMDNLS equation \eqref{foc_CMDNLS} 
with arbitrarily large prescribed finite regularity and decay.
\end{remark}

\begin{proposition}\label{example}
    For every fixed $m\ge1$, there are a nonzero datum
$u_0\in\mathcal S(\mathbb R)\cap H^m_+$ and a time $T_m>0$
such that \eqref{foc_de_CMDNLS} has a nonzero solution
\[
u\in C([-T_m,T_m];H^m_+\cap L_m^2),\qquad u(0)=u_0.
\]
\end{proposition}

\begin{proof}[Proof of Proposition \ref{example}]
Choose a nonzero function $\varphi\in C_c^\infty((1,2))$ and
$\varepsilon>0$, and define $\widehat u_0=\varepsilon\varphi$.
Then $u_0\in\mathcal S(\mathbb R)\cap H^m_+$ and $u_0\ne0$.
The Hardy local theory in \cite[Theorem 2.1]{GerardCPAM} gives
$T_m>0$ and a solution $u\in C([-T_m,T_m];H^m_+)$.
By Theorem~\ref{foc_invariant_HsLr},
$u\in C([-T_m,T_m];H^m_+\cap L_m^2)$.
Since $u(0)=u_0\ne0$, this solution is not identically zero.
\end{proof}

A feature of Theorem~\ref{foc_invariant_HsLr} is that its proof does \emph{not} rely on the complete integrability of the CMDNLS equation; rather, it is based on  energy-type estimates that are applicable to a broader class of differential equations. The strategy proceeds by introducing a family of truncated weight functions $w_N(x)$, which approximate $\langle x\rangle$ while $N\to +\infty$. For these truncated weights, we derive \emph{a priori} estimates on the weighted energy $\|w_N^r u(t)\|_{L^2}^2$ that are uniform in the truncation parameter $N$ by Gr\"{o}nwall's inequality. The passage to the limit $N\to\infty$ is then justified by Fatou's  Lemma, yielding the desired bound for the full weight $\langle x\rangle^r$. 
For $0<r\le1$, the linear term is controlled directly by the
$H^1$ norm. For $r>1$, a weighted interpolation inequality controls
$\|w_N^{r-1}\partial_xu\|_{L^2}$ using the $H^s$ norm and the
truncated weighted norm.

Theorem~\ref{foc_invariant_HsLr} establishes the persistence of the regularity and decay of the solution $u$ itself. However, it does not address the decay of the derivatives $\partial_x^k u$ for $k\ge 1$. Our second result addresses the decay of its derivatives,
using the weighted Sobolev spaces $H_r^s$ defined in \eqref{H_r^s}.

\begin{theorem}\label{H^N_1}
Let $p \ge3$ and $m$ be integers satisfying $1\le m\le p-2$,
and let $0\le r\le p-m$.
Let $0<\mathcal{T}<\mathfrak{T}$.
Suppose that $u\in C^0([-\mathcal{T}, \mathcal{T}],H^{p})$ solves \eqref{foc_de_CMDNLS}
    with initial data $u_0\in H^{p}\cap H_r^m$.
    Then
  \begin{equation}
    u(t)\in H_r^m,\quad \forall t\in [-\mathcal{T}, \mathcal{T}].
    \end{equation} 
Thus $H^p\cap H_r^m$ is invariant under the $H^p$ flow
throughout its lifespan.
\end{theorem}

Theorem~\ref{H^N_1} exploits the Lax pair structure of the CMDNLS equation. 
We characterize $u\in H_r^m$ through the conditions
$\mathbb L_u^k u\in L_r^2$ for every $0\le k\le m$,
under the stated Sobolev regularity assumptions.
The unitary equivalence identifies the evolution of these
iterates, and a separate weighted estimate for that evolution
propagates their spatial decay. The loss of regularity under
iteration explains the condition $r\le p-m$ used in this proof.

It is instructive to compare the two approaches. The proof of Theorem~\ref{foc_invariant_HsLr} uses
weighted energy estimates, interpolation, and Gronwall's
inequality. The proof of Theorem~\ref{H^N_1} combines
weighted energy estimates with the Lax representation
and unitary equivalence.

\subsection{Related work}
De Moura and Pilod proved local well-posedness for the
intermediate nonlinear Schr\"odinger equation and its
deep-water model in $H^s(\mathbb R)$ for $s>1/2$
with arbitrary initial data \cite{deMoura}.
Chapouto, Forlano, and Laurens obtained a full-space local theory
for $s>1/4$ \cite{laurens2025}.
G\'erard and Lenzmann established global well-posedness for the
focusing equation in $H^s_+$ for integer $s\ge1$ when
$\|u_0\|_2^2\le2\pi$, and studied ground states, traveling waves,
and multi-soliton dynamics \cite{GerardCPAM}.
Killip, Laurens, and Vi\c{s}an constructed the flow in $L^2_+$,
with arbitrary mass in the defocusing case and mass below the
soliton mass in the focusing case \cite{KLV2025}.
Hadama proved small-data global well-posedness in the full space
$L^2(\mathbb R)$  \cite{Hadama}.

Kim, Kim, and  Kwon constructed smooth
chiral finite-time blow-up solutions with mass arbitrarily close
to the soliton mass \cite{KimKimKwon}.
Kim and  Kwon established soliton resolution for
finite-time blow-up and global solutions in the setting of
\cite{KimKwon}.
Tian and Zhang studied soliton resolution at the critical
$L^2_+$ regularity for specified classes of initial data and
flow--Lax admissible solutions \cite{Tian}.

For the focusing equation with a nonzero background, Sun
constructed dark multi-soliton solutions and established complete
integrability on the corresponding dark $N$-soliton manifold
\cite{Sun26}.

\subsection{Notation}
\begin{enumerate}
    \item The inverse transform $\check f $ is defined by 
    $\check{f}(x)=\frac{1}{2\pi}\int_\R f(\xi)e^{i\xi x} \d \xi$ for all $f\in L^1$.

    \item We also set $\Pi_-=\mathrm{id}-\Pi$ and
$H^s_-=\{f\in H^s\mid\operatorname{supp}\widehat f
\subset(-\infty,0]\}$.

\item We use $D=-i\partial_x: H^1\to L^2$.

  \item The notation $A\lesssim B$ means that $A\le CB$ for a
positive constant $C$ independent of the truncation parameter $N$.
Any dependence on regularity, weight exponents, or a priori bounds
will be specified in the relevant argument.

  \item For a Banach space $X$, let $\mathbf B_{\mathbb C}(X)$
denote the Banach space of bounded complex-linear operators on $X$,
equipped with the operator norm.
\end{enumerate}

\section{Proof of Theorem \ref{foc_invariant_HsLr}}

\begin{proposition}\label{lem:LWP}
    Let $u_0\in H^s$ and $s>\frac 14$. Then there exists a time 
\begin{equation}\label{frak_T}
\mathfrak{T}=\mathfrak{T}(\| u_0\|_{H^s})>0,
\end{equation}
and a unique solution $u\in C([-\mathfrak{T}, \mathfrak{T}]; H^s)$ of \eqref{foc_de_CMDNLS}  with the initial datum $u_0$.
\end{proposition}
\begin{proof}
    See \cite[Theorem 1.1]{laurens2025} for the details. 
\end{proof}

Now, we give an explicit construction of the truncated weights.
Define
\[
\vartheta(s)=
\begin{cases}
e^{-1/s}, & s>0,\\
0, & s\le0,
\end{cases}
\qquad
\chi(s)=
\frac{\vartheta(2-s)}
{\vartheta(s)+\vartheta(2-s)}.
\]
Then $\chi\in C^\infty(\mathbb R)$ and
\[
0\le\chi\le1,\qquad
\chi(s)=1\quad(s\le0),\qquad
\chi(s)=0\quad(s\ge2).
\]
Moreover, the following identities hold: 
\[
\chi(s)+\chi(2-s)=1,
\qquad
\int_0^2\chi(s)\,\mathrm ds=1.
\]

For any integer $N\ge 1$, set
\[
\delta_N=2N-\langle N\rangle.
\]
Since $\langle N\rangle=\sqrt{1+N^2}$, we have
\[
(2-\sqrt2)N\le\delta_N\le N.
\]
We define
\begin{equation}\label{eq:truncated-weight-definition}
w_N(x)
=
\langle N\rangle
+
\delta_N
\int_0^{(\langle x\rangle-\langle N\rangle)/\delta_N}
\chi(s)\,\mathrm ds.
\end{equation}

If $|x|\le N$, the upper integration limit is nonpositive.
Since $\chi(s)=1$ for $s\le0$, it follows that
\[
w_N(x)
=
\langle N\rangle
+
\delta_N
\frac{\langle x\rangle-\langle N\rangle}{\delta_N}
=
\langle x\rangle.
\]
If $|x|\ge3N$, then
\[
\langle x\rangle
\ge3N
\ge4N-\langle N\rangle
=
\langle N\rangle+2\delta_N.
\]
Thus the upper integration limit is at least $2$, and
\[
w_N(x)
=
\langle N\rangle
+
\delta_N\int_0^2\chi(s)\,\mathrm ds
=
2N.
\]

Differentiating \eqref{eq:truncated-weight-definition}, we obtain
\[
\partial_xw_N(x)
=
\chi\left(
\frac{\langle x\rangle-\langle N\rangle}{\delta_N}
\right)
\frac{x}{\langle x\rangle}.
\]
Hence $w_N$ is smooth and even, is nondecreasing on
$[0,\infty)$, and satisfies $|\partial_xw_N|\le1$.
The inequalities $0\le\chi\le1$ also give
\[
1\le w_N(x)\le\langle x\rangle.
\]
Finally, all positive-order derivatives of $\langle x\rangle$
are bounded. Repeated differentiation therefore yields
\[
\|\partial_x^k w_N\|_{L^\infty}
\le
c_k\sum_{j=0}^{k-1}\delta_N^{-j}
\le C_k,
\qquad k\ge1,
\]
where $c_k$ and $C_k$ depend only on $k$ and the fixed
function $\chi$, and are independent of $N$.
The weights defined above satisfy the truncation properties
in \cite[(2.5)]{Ponce2011}. In particular,
\[
w_N(x)=\langle x\rangle\quad (|x|\le N),\qquad
w_N(x)=2N\quad (|x|\ge3N),
\]
\[
1\le w_N(x)\le2\langle x\rangle,\qquad
|\partial_xw_N(x)|\le1,\qquad
\sup_N\|\partial_x^k w_N\|_{L^\infty}\le C_k\quad(k\ge1).
\]
Here $C_k$ depends only on $k$ and the fixed cutoff function
$\chi$, and is independent of $N$.
Moreover, we have 
\begin{equation}\label{CVPtWise}
    \lim_{N \to + \infty} w_N(x) = \langle x\rangle, \qquad \forall x \in \mathbb{R}.
\end{equation} 
The truncated weighted interpolation inequality in
\cite[Lemma 1]{Ponce2011} holds for this family with a constant
independent of $N$.
For any $a\in\R$, we define the operator  $\mathbb{J}^a = \langle -i \partial_x \rangle^a :H^a\to L^2$ as 
\begin{equation}\label{bfJdef}
   \widehat{\mathbb{J}^a(f)}(\xi)=  (1+\xi^2)^{\frac a 2}\hat{f}(\xi) , \qquad \forall \xi \in \mathbb{R}. 
\end{equation}

\begin{proof}[Proof of Theorem \ref{foc_invariant_HsLr}]
   We only prove the statement for $\sigma=1$. The case $\sigma=-1$ is similar. 
   Throughout the proof, the constants $c_1,c_2,\ldots$ are
positive and independent of $N$ and $t$, but may depend on
$r$, $s$, and the choice of the truncated weights.

We first justify the weighted energy calculation by regularizing
$u$ in time. Choose $\eta\in C_c^\infty((-1,1))$ such that
\(
\eta\ge0 \), \(
\int_{\mathbb R}\eta(\tau)\,\mathrm d\tau=1,
\)
and set
\(
\eta_\varepsilon(\tau)
=
\frac1\varepsilon\eta\left(\frac{\tau}{\varepsilon}\right).
\)
Fix $0<T_0<T$. For $0<\varepsilon<T-T_0$ and
$t\in[-T_0,T_0]$, we define the following Bochner integral: 
\[
u_\varepsilon(t)
=
\int_{-\varepsilon}^{\varepsilon}
\eta_\varepsilon(\tau)u(t-\tau)\,\mathrm d\tau
=
\int_{-T}^{T}\eta_\varepsilon(t-s)u(s)\,\mathrm ds \in H^1.
\]
Differentiating the smooth kernel gives
\[
\partial_tu_\varepsilon(t)
=
\int_{-T}^{T}\eta_\varepsilon'(t-s)u(s)\,\mathrm ds
\in H^1.
\]

We use the duality pairing
$\langle\cdot,\cdot\rangle_{H^{-1},H^1}$ with the convention
\[
\langle F,g\rangle_{H^{-1},H^1}
=
\int_{\mathbb R}F\overline g\,\mathrm dx
\qquad
\text{when }F\in L^2,\quad g\in H^1.
\]
By the Sobolev embedding $H^1\hookrightarrow L^\infty$
and the boundedness of $\mathbb H$ on $L^2$,
\[
\begin{split}
\|u\mathbb H\partial_x(|u|^2)\|_{L^2}
+
\|u\partial_x(|u|^2)\|_{L^2}
\le
2\|u\|_{L^\infty}\|\partial_x(|u|^2)\|_{L^2}
\le
4\|u\|_{L^\infty}^2\|\partial_xu\|_{L^2}.
\end{split}
\]
These nonlinear terms depend continuously on $u\in H^1$.
Since $\partial_x^2:H^1\to H^{-1}$ is bounded, equation \eqref{foc_de_CMDNLS} implies that
\(
\partial_tu\in C([-T,T];H^{-1}).
\)
Consequently, differentiation of the scalar convolution against
any fixed $\phi\in H^1$ yields
\[
\begin{split}
\langle\partial_tu_\varepsilon(t),\phi\rangle_{H^{-1},H^1}
=
\frac{\mathrm d}{\mathrm dt}
\int_{-\varepsilon}^{\varepsilon}
\eta_\varepsilon(\tau)
\langle u(t-\tau),\phi\rangle_{H^{-1},H^1}\,
\mathrm d\tau
=
\int_{-\varepsilon}^{\varepsilon}
\eta_\varepsilon(\tau)
\langle\partial_tu(t-\tau),\phi\rangle_{H^{-1},H^1}\,
\mathrm d\tau.
\end{split}
\]
Thus, we have
\(
\partial_tu_\varepsilon(t)
=
\int_{-\varepsilon}^{\varepsilon}
\eta_\varepsilon(\tau)\partial_tu(t-\tau)\,\mathrm d\tau
\) in $H^{-1}$.
Using $\eta_\varepsilon\ge0$ and
$\int_{\mathbb R}\eta_\varepsilon=1$, we obtain
\[
\begin{split}
\|u_\varepsilon(t)-u(t)\|_{H^1}
&\le
\sup_{|\tau|\le\varepsilon}
\|u(t-\tau)-u(t)\|_{H^1},
\\
\|\partial_tu_\varepsilon(t)-\partial_tu(t)\|_{H^{-1}}
&\le
\sup_{|\tau|\le\varepsilon}
\|\partial_tu(t-\tau)-\partial_tu(t)\|_{H^{-1}}.
\end{split}
\]
The uniform continuity of $u$ in $H^1$ and of $\partial_tu$
in $H^{-1}$ therefore implies
\[
\begin{split}
u_\varepsilon&\longrightarrow u
\quad\text{in }C([-T_0,T_0];H^1),
\\
\partial_tu_\varepsilon&\longrightarrow\partial_tu
\quad\text{in }C([-T_0,T_0];H^{-1}).
\end{split}
\]

For each fixed $N$, the functions $w_N^{2r}$ and
$\partial_x(w_N^{2r})$ are bounded. The product rule implies that
\[
\partial_x(w_N^{2r}f)
=
w_N^{2r}\partial_xf+\partial_x(w_N^{2r})f.
\]
Then there exists a constant $c_{1,N} > 0$,  which  depends on $N$ and $r$, such that
\[
\|w_N^{2r}f\|_{H^1}
\le c_{1,N}\|f\|_{H^1},
\qquad \forall f\in H^1.
\]
 
Since $u_\varepsilon$ is smooth in time with values in $H^1$,
the usual differentiation rule gives
\[
\begin{split}
\frac12\frac{\mathrm d}{\mathrm dt}
\|w_N^ru_\varepsilon(t)\|_{L^2}^2
=
\operatorname{Re}\int_{\mathbb R}
\partial_tu_\varepsilon(t)\,
w_N^{2r}\overline{u_\varepsilon(t)}\,\mathrm dx
=
\operatorname{Re}
\langle
\partial_tu_\varepsilon(t),w_N^{2r}u_\varepsilon(t)
\rangle_{H^{-1},H^1}.
\end{split}
\]
Integrating from $0$ to $t\in[-T_0,T_0]$, we obtain
\[
\begin{split}
&\frac12\|w_N^ru_\varepsilon(t)\|_{L^2}^2
-\frac12\|w_N^ru_\varepsilon(0)\|_{L^2}^2
=
\operatorname{Re}\int_0^t
\langle
\partial_\tau u_\varepsilon(\tau),
w_N^{2r}u_\varepsilon(\tau)
\rangle_{H^{-1},H^1}\,\mathrm d\tau.
\end{split}
\]
The convergence of the integrand follows from
\[
\begin{split}
&
\left|
\langle\partial_tu_\varepsilon,w_N^{2r}u_\varepsilon\rangle_{H^{-1},H^1}
-
\langle\partial_tu,w_N^{2r}u\rangle_{H^{-1},H^1}
\right|
\le
c_{1,N}
\|\partial_tu_\varepsilon-\partial_tu\|_{H^{-1}}
\|u_\varepsilon\|_{H^1}
+
c_{1,N}
\|\partial_tu\|_{H^{-1}}
\|u_\varepsilon-u\|_{H^1}
\longrightarrow0,
\end{split}
\]
uniformly for $t\in[-T_0,T_0]$.
Moreover, boundedness of $w_N^r$ gives
\[
\|w_N^r(u_\varepsilon-u)\|_{L^2}
\le
\|w_N^r\|_{L^\infty}\|u_\varepsilon-u\|_{L^2}
\longrightarrow0
\]
uniformly on the same interval.
Passing to the limit $\varepsilon\to0$, we conclude that
\[
\begin{split}
&\frac12\|w_N^ru(t)\|_{L^2}^2
-\frac12\|w_N^ru(0)\|_{L^2}^2
=
\operatorname{Re}\int_0^t
\langle
\partial_\tau u(\tau),w_N^{2r}u(\tau)
\rangle_{H^{-1},H^1}\,\mathrm d\tau.
\end{split}
\]
The integrand is continuous in time.
Since $T_0<T$ is arbitrary, the fundamental theorem of calculus
yields 
that the function
$t \in (-T,T) \mapsto \|w_N^ru(t)\|_{L^2}^2$ is differentiable and 
\begin{equation}\label{ddtL2wNrut}
\frac12\frac{\mathrm d}{\mathrm dt}
\|w_N^ru(t)\|_{L^2}^2
=
\operatorname{Re}
\langle
\partial_tu(t),w_N^{2r}u(t)
\rangle_{H^{-1},H^1},
\qquad t\in(-T,T).
\end{equation}
The integrated identity extends to $t=\pm T$ by continuity.

We next justify the spatial integration by parts.
By the definition of the distributional derivative and the
density of $C_c^\infty(\mathbb R)$ in $H^1$,
\[
\langle\partial_x^2u,\phi\rangle_{H^{-1},H^1}
=
-\int_{\mathbb R}\partial_xu\,
\overline{\partial_x\phi}\,\mathrm dx,
\qquad \phi\in H^1.
\]
Taking $\phi=w_N^{2r}u\in H^1$, we obtain
\begin{equation}\label{Hneg1H1LapuwN2r}
\begin{split}
\langle
\partial_x^2u,w_N^{2r}u
\rangle_{H^{-1},H^1}
=
-\int_{\mathbb R}
\partial_xu\,
\partial_x\bigl(w_N^{2r}\overline u\bigr)\,
\mathrm dx
=
-\int_{\mathbb R}w_N^{2r}|\partial_xu|^2\,\mathrm dx
-
\int_{\mathbb R}
\partial_x(w_N^{2r})\overline u\,\partial_xu\,
\mathrm dx.
\end{split}
\end{equation}
Both integrals are finite, since
\[
\begin{split}
\int_{\mathbb R}w_N^{2r}|\partial_xu|^2\,\mathrm dx
&\le
\|w_N^{2r}\|_{L^\infty}\|\partial_xu\|_{L^2}^2,
\\
\left|
\int_{\mathbb R}
\partial_x(w_N^{2r})\overline u\,\partial_xu\,
\mathrm dx
\right|
&\le
\|\partial_x(w_N^{2r})\|_{L^\infty}
\|u\|_{L^2}\|\partial_xu\|_{L^2}.
\end{split}
\]

Finally, substituting \eqref{foc_CMDNLS} into
\eqref{ddtL2wNrut}, with the linear term interpreted through
\eqref{Hneg1H1LapuwN2r}, gives
\begin{equation}\label{eq:multi_w_2r_u}
\begin{split}
& \frac12\frac{\mathrm d}{\mathrm dt}\|w_N^ru\|_{L^2}^2 \\
 = &
\operatorname{Re}\left(
i\langle
\partial_x^2u,w_N^{2r}u
\rangle_{H^{-1},H^1}
\right)
+
\operatorname{Re}\left(
i\int_{\mathbb R}
w_N^{2r}|u|^2
\mathbb H\bigl(\partial_x|u|^2\bigr)\,
\mathrm dx
\right)
+
\operatorname{Re}\left(
\int_{\mathbb R}
w_N^{2r}|u|^2\partial_x(|u|^2)\,
\mathrm dx
\right).
\end{split}
\end{equation}
Since $\mathbb{H}(\partial_x|u|^2)$ is real-valued,
the second term on the right-hand side of
\eqref{eq:multi_w_2r_u} vanishes.
For the last term, integration by parts gives
\begin{equation}
    \Re\left(
        \int_{\mathbb{R}}
        w_N^{2r}|u|^2\partial_x(|u|^2)\,\mathrm{d}x
    \right)
    =
    -\int_\R \partial_x(w_N^{2r})|u|^4 \d x -\int_\R w_N^{2r} |u|^2 \partial_x(|u|^2) \d x, 
\end{equation}
This implies that 
\begin{equation}
\begin{split}
    \Re\left(
        \int_{\mathbb{R}}
        w_N^{2r}|u|^2\partial_x(|u|^2)\,\mathrm{d}x
    \right)
    =
    -\frac12\int_{\mathbb{R}}
        \partial_x(w_N^{2r})|u|^4\,\mathrm{d}x
    =
    -r\int_{\mathbb{R}}
        w_N^{2r-1}(\partial_xw_N)|u|^4\,\mathrm{d}x.
\end{split}
\end{equation}
It follows from $w_N\ge1$, 
$\partial_xw_N \le 1$, and the Sobolev embedding
$H^s(\mathbb{R})\hookrightarrow L^\infty(\mathbb{R})$
that
\begin{equation}\label{last_term}
\begin{split}
    \left|
    \Re\left(
        \int_{\mathbb{R}}
        w_N^{2r}|u|^2\partial_x(|u|^2)\,\mathrm{d}x
    \right)
    \right|
    \le
    c_1\|u\|_{L^\infty}^2\|w_N^r u\|_{L^2}^2
    \le
    c_2\|u\|_{H^s}^2\|w_N^r u\|_{L^2}^2.
\end{split}
\end{equation}

We next estimate the first term on the right-hand side of
\eqref{eq:multi_w_2r_u}. Integrating by parts,
we have
\begin{equation}
\begin{split}
    I
    &:=
    \Re \big( i \langle
\partial_x^2u,w_N^{2r}u
\rangle_{H^{-1},H^1}
\big)
    =
    \Im\left(
        \int_{\mathbb{R}}
        \partial_x(w_N^{2r})\overline{u}\,
        \partial_xu\,\mathrm{d}x
        +
        \int_{\mathbb{R}}
        w_N^{2r}|\partial_xu|^2\,\mathrm{d}x
    \right)
    \\
    &=
    2r\Im\left(
        \int_{\mathbb{R}}
        w_N^{2r-1}(\partial_xw_N)
        \overline{u}\,\partial_xu\,\mathrm{d}x
    \right).
\end{split}
\end{equation}
Thus, by the Cauchy--Schwarz inequality,
\begin{equation}\label{I_esti}
    |I|
    \le
    c_3\|w_N^r u\|_{L^2}
        \|w_N^{r-1}\partial_xu\|_{L^2}.
\end{equation}

For $0<r\le1$, the estimate \eqref{I_esti} and $w_N\ge1$ give
\[
|I|\le c_3\|w_N^ru\|_2\|\partial_xu\|_2.
\]
We use weighted interpolation only in the remaining case $r>1$.
The Leibniz formula yields
\begin{equation}\label{w_r-1_partial_u}
\begin{split}
    \|w_N^{r-1}\partial_xu\|_{L^2}
    \le
    \|\partial_x(w_N^{r-1}u)\|_{L^2}
    +
    \|(\partial_xw_N^{r-1})u\|_{L^2}
    \le
    \|\mathbb{J}(w_N^{r-1}u)\|_{L^2}
    +
    \|(\partial_xw_N^{r-1})u\|_{L^2}.
\end{split}
\end{equation}
Applying the weighted interpolation inequality (see Lemma 1 in \cite{Ponce2011}) with
$a=b=r$ and $\theta=1/r$, we obtain
\begin{equation}\label{J_w_r-1_u}
\begin{split}
    \|\mathbb{J}(w_N^{r-1}u)\|_{L^2}
    \le
    c_4\|\mathbb{J}^r u\|_{L^2}^{1/r}
        \|w_N^r u\|_{L^2}^{(r-1)/r}
    \le
    c_4\|u\|_{H^s}^{1/r}
        \|w_N^r u\|_{L^2}^{(r-1)/r},
\end{split}
\end{equation}
where the last inequality follows from $r\le s$. 
Moreover,
\begin{equation}
    |\partial_x(w_N^{r-1})|
    =
    (r-1)w_N^{r-2}|\partial_xw_N|
    \le c_5 w_N^{r-2}.
\end{equation}
Since $w_N\ge1$, we have $w_N^{r-2}\le w_N^r$, and hence
\begin{equation}\label{par_x_r-1_u}
    \|(\partial_xw_N^{r-1})u\|_{L^2}
    \le
    c_5\|w_N^{r-2}u\|_{L^2}
    \le
    c_5\|w_N^r u\|_{L^2}.
\end{equation}
Combining \eqref{w_r-1_partial_u} \eqref{J_w_r-1_u} \eqref{par_x_r-1_u}, we conclude that
\begin{equation}
    \|w_N^{r-1}\partial_xu\|_{L^2}
    \le
    c_4\|u\|_{H^s}^{1/r}
        \|w_N^r u\|_{L^2}^{(r-1)/r}
    +
    c_5\|w_N^r u\|_{L^2}.
\end{equation}
Substituting this estimate into
\eqref{I_esti}, we obtain
\begin{equation}\label{I_esti_2}
    |I|
    \le
    c_6\|u\|_{H^s}^{1/r}
        \|w_N^r u\|_{L^2}^{2-1/r}
    +
    c_7\|w_N^r u\|_{L^2}^2,
\end{equation}
where $c_6=c_3c_4$ and $c_7=c_3c_5$.
By Young's inequality with conjugate exponents
$2r$ and $2r/(2r-1)$,
\begin{equation}\label{young}
    c_6\|u\|_{H^s}^{1/r}
        \|w_N^r u\|_{L^2}^{2-1/r}
    \le
    c_8\|u\|_{H^s}^2
    +
    c_9\|w_N^r u\|_{L^2}^2.
\end{equation}
Therefore, \eqref{eq:multi_w_2r_u},
\eqref{I_esti_2}, \eqref{young} and
\eqref{last_term} imply that
\begin{equation}\label{s>1_esti_zong}
\begin{split}
    \left|
        \partial_t\|w_N^r u(t)\|_{L^2}^2
    \right|
    &\le
    2c_8\|u(t)\|_{H^s}^2
    +
    \left(
        2c_7+2c_9+2c_2\|u(t)\|_{H^s}^2
    \right)
    \|w_N^r u(t)\|_{L^2}^2.
\end{split}
\end{equation}

By Proposition~\ref{lem:LWP}, we may define
\[
    M_T:=\sup_{-T\le t\le T}\|u(t)\|_{H^s}<+\infty.
\]
Set
\[
    U(t):=\|w_N^r u(t)\|_{L^2}^2,
    \qquad
    A_T:=2c_7+2c_9+2c_2M_T^2,
    \qquad
    B_T:=2c_8M_T^2.
\]
Then \eqref{s>1_esti_zong} gives
\begin{equation}
    |\partial_tU(t)|
    \le A_TU(t)+B_T.
\end{equation}
For $t\in[0,T]$, multiplying by $e^{-tA_T}$ and
integrating from $0$ to $t$, Gr\"onwall's inequality gives that 
\begin{equation}
\begin{split}
    U(t)
    &\le
    e^{tA_T}U(0)
    +
    B_T\int_0^t e^{(t-\tau)A_T}\,\mathrm{d}\tau
    \le
    e^{tA_T}\bigl(U(0)+B_T(-\frac{1}{A_T})(e^{-tA_T}-1)\bigr)
    \le 
    e^{tA_T}\bigl(U(0)+\frac{B_T}{A_T} \bigr)
    .
\end{split}
\end{equation}
For $t\in[-T,0]$, applying the similar argument to
$U(-t)$, we have 
\begin{equation}
    U(t)
    \le
    e^{-tA_T}\bigl(U(0)+\frac{B_T}{A_T}\bigr).
\end{equation}
Consequently, it follows that 
\begin{equation}\label{uniform_w_r_u}
    \sup_{-T\le t\le T}\|w_N^r u(t)\|_{L^2}^2
    \le
    e^{TA_T}
    \left(
        \|w_N^r u_0\|_{L^2}^2+\frac{B_T}{A_T}
    \right).
\end{equation}

The truncated weights satisfy
$w_N(x)\le2\langle x\rangle$ and
$w_N(x)\to\langle x\rangle$ as $N\to\infty$.
Thus,
\[
    \|w_N^r u_0\|_{L^2}^2
    \le
    c_{10}\|\langle x\rangle^r u_0\|_{L^2}^2,
    \qquad c_{10}=2^{2r}.
\]
For each $t\in[-T,T]$, Fatou's lemma and
\eqref{uniform_w_r_u} now yield
\begin{equation}
\begin{split}
    \|\langle x\rangle^r u(t)\|_{L^2}^2
    \le
    \liminf_{N\to\infty}\|w_N^r u(t)\|_{L^2}^2
    \le
    e^{TA_T}
    \left(
        c_{10}\|\langle x\rangle^r u_0\|_{L^2}^2
        +\frac{B_T}{A_T}
    \right)
    <+\infty.
\end{split}
\end{equation}
It follows that $u(t)\in L_r^2$ for all $t\in[-T,T]$.
Together with $u\in C([-T,T];H^s)$, this proves that
\[
    u(t)\in H^s\cap L_r^2,
    \qquad \forall\,t\in[-T,T].
\]
The energy bound  implies
$|U(t)-U(t_0)|\le c_T|t-t_0|$, uniformly in $N$.
Dominated convergence at each fixed time shows that
$t\mapsto\|\langle x\rangle^r u(t)\|_2^2$ is continuous.

\end{proof}

\section{Non-chiral Lax pair and unitary propagator}

Let $u\in H^1$. Define the Lax pair on the full Sobolev space $H^1$ as 
     \begin{align}
    &
        \mathcal{L}_u:= D-\sigma m_u\mathcal{T}_{\overline u}: H^1\to L^2,
    \\
    &
        \mathcal{B}_u:=\sigma m_u\mathcal{T}_{\partial_x\overline u}
        -
       \sigma  m_{\partial_x u}\mathcal{T}_{\overline u}
        +
        i(m_u \mathcal{T}_{\overline u})^2: H^1\to L^2,
    \end{align}
where $m_u: L^2\to L^2$ denotes the multiplication operator and for any $h\in L^2$, $m_u(h):=uh$, and define $\mathcal{T}_u:=\Pi\circ m_u: L^2\to L^2_+$ as $\mathcal{T}_u(h):=\Pi(uh)$. If $u$ solves \eqref{foc_de_CMDNLS}, then we have 
\begin{equation}
    \partial_t \mathcal{L}_u=[\mathcal{B}_u, \mathcal{L}_u].
\end{equation}
See the proof in \cite[Proposition 2.1]{Ba2024}.

Let $u\in C([-T,T]; H^2)$. Since $\mathcal{B}_u$ is a bounded and skew-adjoint operator on $L^2$, i.e., $\mathcal{B}_u^*=-\mathcal{B}_u$. 
Let $W$ be the unique solution of
\[
\partial_tW(t)=\mathcal B_{u(t)}W(t),\qquad W(0)=\mathrm{id}.
\]
Then 
\begin{equation}
    W\in  C^1\big( 
    [-T, T]; \mathbf{B}_\C(L^2)
    \big)
\end{equation}
satisfying the following properties 
\begin{enumerate}
    \item For any $t\in (-T, T)$, $W(t)$ is a unitary operator on $L^2$, i.e.,
    \begin{equation}
        W(t)W(t)^*=W(t)^*W(t)=\mathrm{id}_{L^2},
        \quad 
        W(t)^*=W(t)^{-1},
    \end{equation}
where $W(t)^*$ denotes the $L^2$-adjoint of $W(t)\in \mathbf{B}_\C(L^2)$.

\item For any $t\in (-T, T)$, we have 
\begin{equation}
    \begin{split}
    &
        \lim_{r\to 0} \| \frac{W(t+r)-W(t)}{r} -\mathcal{B}_{u(t)} W(t)  \|_{\mathbf{B}_\C(L^2)}=0,
        \\
        &
        \lim_{r\to 0} \| \frac{W(t+r)^{-1}- W(t)^{-1}}{r}+W(t)^{-1} \mathcal{B}_{u(t)} \|_{\mathbf{B}_\C(L^2)}=0.
    \end{split}
\end{equation}

\end{enumerate}
For any  $t\in[-T, T]$, we have the unitary equivalence
\begin{equation}
    \mathcal{L}_{u(t)} W(t)=W(t) \mathcal{L}_{u(0)}: H^1 \to L^2.
\end{equation}
See details in \cite[Sec. 4.2]{Sun26}.

\section{Proof of Theorem \ref{H^N_1}}

\begin{lemma}\label{Hns_Hn+1_inHns}
    For every integer $n\ge0$ and $r\ge 0$, we have
    \begin{equation}
    H^n_r\cdot H^{n+1}\subset H^n_r.
    \end{equation}
\end{lemma}
\begin{proof}
    We first verify the case $n=1$. For any $f\in H^1_r$ and $g\in H^2$, we have 
    \begin{equation}
        fg\in L^2_r\cdot L^\infty \subset L^2_r,\quad \partial_x(fg)=f_xg+fg_x\in L^2_r\cdot L^\infty+L^2_r\cdot L^\infty\subset L^2_r.
    \end{equation}
    Hence  $H^1_r\cdot H^2\subset H^1_r$. 
   Now suppose that $H^k_r\cdot H^{k+1}\subset H^k_r$ holds for some $k\in\mathbb{N}_+$ and all $r>0$. 
   For any $f\in H^{k+1}_r$ and $g\in H^{k+2}$, it follows that 
   \begin{equation}
       fg\in H^k_r\cdot H^{k+1}\subset H^k_r,\quad \partial_x(fg)=f_xg+fg_x\in H^k_r\cdot H^{k+1}+H^k_r\cdot H^{k+1}\subset H^k_r.
   \end{equation}
Thus $fg\in H^{k+1}_r$, which yields $H^{k+1}_r\cdot H^{k+2}\subset H^{k+1}_r$. By mathematical induction,
   the inclusion $H^n_r\cdot H^{n+1}\subset H^n_r$ is valid for every $n\in\mathbb{N}_+$ and every $r>0$.
   
\end{proof}

\begin{lemma}\label{LuH_m+1_in_Hm}
Let $p\ge 1$  be an integer and $r>0$. Assume that $u\in H^p_r$. Then 
    \begin{equation}
        \mathcal{L}_u(H^{m+1}_r)\subset H^m_r, \qquad \text{for every integer }  m \text{ with } 0 \leq m \leq  p-1.
        \end{equation}
\end{lemma}

\begin{proof}
    Fix  an integer $m$ with $0\leq m \leq p-1$. For any $f\in H^{m+1}_r$, we  have $\overline{u} f \in H^n_r\cdot H^{m+1}_r \subset H^{m+1}$, which implies $\Pi(\overline u f)\in H^{m+1}_+$. By Lemma~\ref{Hns_Hn+1_inHns}, it follows that $u\Pi(\overline u f)\in H^p_r\cdot H^{m+1}\subset H^m_r \cdot H^{m+1}\subset H^m_r$. Therefore, we obtain 
    \begin{equation}
        \mathcal{L}_u(f)= -i\partial_x f-\sigma u\Pi(\overline{u} f)\in H^m_r+H^m_r\subset H^m_r.
    \end{equation}
    Hence $\mathcal{L}_u (H^{m+1}_r)\subset H^m_r$ for all integers $m$ with $0\le m\le p-1$.
\end{proof}

\begin{proposition}\label{pro-u-to-Lu}
    Let $p\in \mathbb{N}_+$ and assume that $u\in H^p_r$ for any $r>0$. Then 
    \begin{equation}\label{u-to-Lu}
        \mathcal{L}_u(u),\, \mathcal{L}^2_u(u),\, \dots, \mathcal{L}_u^p(u)\in L^2_r.
    \end{equation} 
\end{proposition}

\begin{proof}
    For every integer $k$ with $0\leq k\leq p$, we claim that 
    \begin{equation}\label{claim}
        v_{k}:= \mathcal{L}_u^{k}(u)\in H^{p-k}_r.
    \end{equation}
We prove this by induction on $k$. For $k=0$, the claim is trivially true because $v_0=u\in H^p_r$ by hypothesis. 
    For $k=1$, we have $v_1=\mathcal{L}_u(u)\in \mathcal{L}_u(H^p_r)\subset H^{p-1}_r$ by Lemma~\ref{LuH_m+1_in_Hm}.  Now assume that \eqref{claim} holds for some  integer $k$ with $1 \le k\leq p-1$. Then we obtain 
    \begin{equation}
    v_{k+1}=\mathcal{L}_u^{k+1}(u)=\mathcal{L}_u(v_{k})\in \mathcal{L}_u(H^{p-k}_r)\subset H^{p-k-1}_r,
    \end{equation}
    which is exactly the desired statement for $k+1$. 
    Thus the claim is valid for $0\le k \le  p$.  In particular, for $k=1,2,\dots, p$, we have $v_k\in H^{p-k}_r$, i.e., $\mathcal{L}_u(u), \mathcal{L}^2_u(u), \cdots, \mathcal{L}_u^p(u)\in L^2_r$.
\end{proof}

\begin{lemma}\label{Lu(fg)=}
    If $u, f, g\in H^1$, then 
    \begin{equation}
        \mathcal{L}_u(fg)=(\mathcal{L}_u f) g+f(Dg)-\sigma u\Pi\bigl(g\Pi_-(\overline{u} f)\bigr)+\sigma u \Pi_- (g \Pi(\overline u f)),
    \end{equation}
    where $\Pi_-=\mathrm{id}-\Pi: L^2\to L^2_-$.
\end{lemma}

\begin{proof} Since $g\in H^1$, a direct computation shows that
    \begin{equation}
    \begin{split}
        \mathcal{L}_u(fg)
        &
        =
        D(fg)-\sigma u\Pi(\overline{u}fg) =(Df)g+f(Dg)-\sigma u\Pi(\overline{u} fg)
        \\
        &=\bigl(Df-\sigma u\Pi(\overline{u}f)\bigr)g+f(Dg)-\sigma u\Pi(\overline{u}fg)+\sigma ug\Pi(\overline{u}f).
    \end{split}
    \end{equation}
Since $g\Pi(\overline u  f)= \Pi(g\Pi(\overline u f))+\Pi_-(g\Pi(\overline u f ))$, 
we have 
\begin{equation}
\begin{split}
\mathcal{L}_u(fg)
&
        =
        (\mathcal{L}_uf)g+f(Dg)-\sigma u\Pi\bigl(\overline{u}fg-g\Pi(\overline u f)\bigr)+\sigma u \Pi_- (g \Pi(\overline u f)) 
        \\
        &
        =
        (\mathcal{L}_u f)g+f(Dg)-\sigma u\Pi\bigl(g\, \Pi_-(\overline{u}f)\bigr)+\sigma u \Pi_- (g \Pi(\overline u f)) .
\end{split}
\end{equation}
\end{proof}

Let $p\in \mathbb{N}_+$ and $r>0$. Assume that $u\in H^p$ such that $u,\, \mathcal{L}_u(u), \, \mathcal{L}_u^2(u), \, \dots, \mathcal{L}_u^p(u)\in L_r^2$. For an integer $k$ with $0\le k\le p$, define a space 
\begin{equation}
    Y^k_{r,u}=\{ f\in H^k: f, \mathcal{L}_u(f), \mathcal{L}_u^2(f), \dots, \mathcal{L}_u^k(f)\in L^2_r \}.
\end{equation}

\begin{lemma}\label{YkHk+1_in_Yk} Assume that $u\in H^p$ such that $u,\, \mathcal{L}_u(u), \, \mathcal{L}_u^2(u), \, \dots, \mathcal{L}_u^p(u)\in L_r^2$.
   For every integer $k$ with $0\le k\le p$,  we have
    \begin{equation}
        Y^k_{r,u}\cdot H^{k+1}\subset Y^k_{r,u}.
    \end{equation}
\end{lemma}

\begin{proof}
    We prove the lemma by induction on $k$. For $k=0$, we have 
    \begin{equation}
        \Y^0\cdot H^1\subset (L^2\cdot H^1)\cap (L^2_r\cdot L^\infty)\subset  L^2 \cap  L^2_r= \Y^0. 
    \end{equation}
    For $k=1$, take any $f\in \Y^1$ and  $g\in H^2$. Then   $fg\in (H^1\cdot H^2)\cap (L^2_r\cdot H^2)\subset H^1\cap L^2_r$. Lemma~\ref{Lu(fg)=} implies that
    \begin{equation}
    \begin{split}
               \mathcal{L}_u(fg)&=\mathcal{L}_u(f)g+f(Dg)-\sigma u\Pi\bigl(g\Pi_-(\overline{u}f)\bigr) +\sigma u \Pi_- (g \Pi(\overline u f))
               \\
               &
               \in L^2_r\cdot H^2+L^2_r\cdot H^1 +L^2_r\cdot \Pi\bigl(H^2\cdot \Pi_-(H^p\cdot H^1)\bigr)
               +L_r^2 \cdot \Pi_- (H^1 \cdot \Pi(H^p \cdot H^1))
               \\
               &
               \subset
               L_r^2+L_r^2+L_r^2\cdot H^1 + L_r^2 \cdot H^1
               \subset L^2_r.
    \end{split}
    \end{equation}
    Hence $fg\in \Y^1$. Then $\Y^1\cdot H^2\subset \Y^1$.
Now  assume that $\Y^n\cdot H^{n+1} \subset \Y^n$ holds for  integers $n$ with $0\le n \le p-1$. Let $f\in \Y^{n+1}$ and $g\in H^{n+2}$. As before, $fg\in (H^{n+1}\cdot H^{n+2}) \cap (L_r^2\cdot H^{n+2})\subset H^{n+1} \cap L^2_r$. We claim that 
\begin{equation}
    \mathcal{L}_u(fg)\in \Y^n. 
\end{equation}
\begin{proof}[Proof of the claim]
    Observe that $\mathcal{L}_u(f)\in \mathcal{L}_u(\Y^{n+1})\subset \Y^n$. By hypothesis, we have 
    \begin{equation}
        (\mathcal{L}_u f)g\in \Y^{n}\cdot H^{n+2}\subset \Y^n\cdot H^{n+1}\subset \Y^n,
    \end{equation}
    and 
    \begin{equation}
        f(Dg)\in \Y^{n+1}\cdot H^{n+1}\subset \Y^n \cdot H^{n+1} \subset \Y^n.
    \end{equation}
    Since $\Pi(g\, \Pi_-(\overline{u}f))\in \Pi(H^{n+2}\cdot \Pi_-\bigl(H^p\cdot H^{n+1})\bigr)\subset \Pi(H^{n+2}\cdot H^{n+1})\subset H^{n+1}$, it follows that
    \begin{equation}
        u\Pi(g\, \Pi_-(\overline{u}f))\in \Y^p\cdot H^{n+1}\subset \Y^n \cdot H^{n+1}\subset \Y^n.
    \end{equation}
    Similarly, since $n+1\le p$, the boundedness of $\Pi$
and $\Pi_-$ on $H^{n+1}$ gives
\[
\Pi_-\bigl(g\Pi(\overline u f)\bigr)\in H^{n+1}.
\]
Using $u\in Y_{r,u}^p\subset Y_{r,u}^n$ and the
induction hypothesis, we obtain
\[
u\Pi_-\bigl(g\Pi(\overline u f)\bigr)
\in Y_{r,u}^n\cdot H^{n+1}\subset Y_{r,u}^n.
\]
    Subsequently, we obtain 
    \begin{equation}
        \mathcal{L}_u(fg)=\mathcal{L}_u(f)g + f(Dg) -\sigma u\Pi(g\, \Pi_-(\overline{u}f))+\sigma u \Pi_- (g \Pi(\overline u f)) \in \Y^n.
    \end{equation}
    The claim is obtained.
\end{proof}
The claim implies $\mathcal{L}_u(fg), \mathcal{L}_u^2(fg), \dots , \mathcal{L}_u^{n+1}(fg)\in L^2_r$. We conclude $fg\in \Y^{n+1}$. Finally, Lemma~\ref{YkHk+1_in_Yk} is obtained by induction.
\end{proof}

\begin{lemma}
Assume that $u\in H^p$ such that $u,\, \mathcal{L}_u(u), \, \mathcal{L}_u^2(u), \, \dots, \mathcal{L}_u^p(u)\in L_r^2$.
   For every integer $k$ with $0\le k\le p-1$, we have
    \begin{equation}
        D\Y^{k+1}\subset \Y^k.
    \end{equation}
\end{lemma}

\begin{proof}
    Let $\phi\in \Y^{k+1}$, we have $D\phi=\mathcal{L}_u(\phi)+\sigma u\Pi(\overline{u} \phi)\in H^k$. By Lemma~\ref{YkHk+1_in_Yk}, it follows that
    \begin{equation}
        u\Pi(\overline{u}\phi)\in \Y^p\cdot \Pi(H^p\cdot H^{k+1})\subset \Y^p\cdot H^{k+1} \subset \Y^k\cdot H^{k+1}\subset \Y^{k}.
    \end{equation}
    Since $\mathcal{L}_u(\phi)\in \mathcal{L}_u(\Y^{k+1})\subset \Y^k$, we conclude that $D\phi= \mathcal{L}_u(\phi)+\sigma u\Pi(\overline{u}\phi)\in \Y^k$.
\end{proof}

\begin{proposition}\label{prop_Lu-to-u}
    Let $p\in \mathbb{N}_+$ and $r>0$. Assume that $u\in H^p$ such that $u, \mathcal{L}_u(u), \mathcal{L}^2_u(u), \dots, \mathcal{L}_u^p(u)\in L_r^2$. Then $u\in H^p_r$.
\end{proposition}

\begin{proof}
    By mathematical induction, for all $m\in [1,p]\cap \mathbb{Z}$, we have 
    \begin{equation}
        D^m(\Y^p)\subset \Y^{p-m}\subset L_r^2.
    \end{equation}
    In particular, $D^mu\in D^m(\Y^p)\subset L_r^2$ for every $m\in [1,p]\cap \mathbb{Z}$. Subsequently, $u\in H^p_r$.
\end{proof}

\begin{lemma}\label{W_exp_Lr2}
      For  $p\ge  2$ and $0\le r\le p$. Let $u\in C^0\big( [-T, T]; H^p \big)$ solve \eqref{foc_de_CMDNLS} with $u(0)=u_0\in H^p\cap L_r^2$. And for $0<\mathcal T\le T$,   let $\psi\in C^0\big( [-\mathcal T, \mathcal T]; H^p  \big)$ solve 
    \begin{equation}\label{eq:psi}
        \partial_t \psi = i\partial_x^2\psi + 2 \sigma u \partial_x \Pi(\overline u \psi), \quad (t,x)\in [-\mathcal T, \mathcal T]\times \R,
    \end{equation}
    with initial datum $\psi(0,x)=\psi_0\in H^p$.
    
    If $\psi_0\in L_r^2$,  then  there exists 
    \begin{equation}
        \psi(t)\in L_r^2, \quad  t\in [-\mathcal T, \mathcal T
        ].
    \end{equation}
\end{lemma}

\begin{proof}
We give the argument for $\sigma=1$, and the case $\sigma=-1$ is
identical. 
Multiplying $w_N^{2r}\overline{\psi}$ on the both two sides of \eqref{eq:psi}, integrating on the real line, 
 and then taking real parts, we obtain
\begin{equation}\label{eq:weighted-energy-identity}
  \frac12\frac{d}{dt}\|w_N^r\psi(t)\|_{L^2}^2
  =I_N(t)+II_N(t),
\end{equation}
where
\begin{align}
  I_N(t)
  :=\operatorname{Re}\int_{\mathbb R}
       i\,\partial_x^2\psi\,w_N^{2r}\overline{\psi}\,\d x,
      \quad 
  II_N(t)
  :=2\operatorname{Re}\int_{\mathbb R}
       u\,\partial_x\Pi(\overline{u}\psi)
       \,w_N^{2r}\overline{\psi}\,\d x.
\end{align}
Here and below, all functions in the integrals are evaluated at time
$t$. The constants $d_1, d_2, d_3\ldots$ are independent on $N$.
We begin with $I_N$. An integration by parts gives
\begin{align}
 I_N(t)
 = -\operatorname{Im}\int_{\mathbb R}
        \partial_x^2\psi\,w_N^{2r}\overline{\psi}\,\d x
 =\operatorname{Im}\int_{\mathbb R}
        \partial_x\psi\,\partial_x(w_N^{2r}\overline{\psi})\,\d x
       \notag
 =\operatorname{Im}\int_{\mathbb R}
        \partial_x\psi\,\partial_x(w_N^{2r})
        \,\overline{\psi}\,\d x,
        \label{eq:linear-integration-by-parts}
\end{align}
since
$\operatorname{Im}\int_{\mathbb R}w_N^{2r}|\partial_x\psi|^2\,\d x=0$.
Using
the Cauchy--Schwarz inequality and
$|\partial_xw_N(x)|\le 1$, we deduce that
\begin{equation}\label{eq:linear-initial-bound}
  |I_N(t)|
  \leq d_1 \|w_N^r\psi(t)\|_{L^2}
           \|w_N^{r-1}\partial_x\psi(t)\|_{L^2}.
\end{equation}
We estimate the last factor by considering two cases.

\smallskip
\noindent\emph{Case 1: $0<r\leq1$.}
Since $w_N\geq1$, we have $w_N^{r-1}\leq1$, and hence
\begin{equation}\label{eq:linear-case-small-r}
  |I_N(t)|
  \leq d_1\|w_N^r\psi(t)\|_{L^2}
           \|\partial_x\psi(t)\|_{L^2}
  \leq d_1\|w_N^r\psi(t)\|_{L^2}
           \|\psi(t)\|_{H^1}.
\end{equation}

\smallskip
\noindent\emph{Case 2: $r>1$.}
The product rule yields
\begin{equation}\label{eq:weighted-product-rule}
  w_N^{r-1}\partial_x\psi
   =\partial_x(w_N^{r-1}\psi)
    -(r-1)w_N^{r-2}(\partial_xw_N)\psi.
\end{equation}
Consequently, by the triangle inequality and the bound on
$\partial_xw_N$,
\begin{align}
 \|w_N^{r-1}\partial_x\psi\|_{L^2}
 \leq\|\partial_x(w_N^{r-1}\psi)\|_{L^2}
      +d_2 \|w_N^{r-2}\psi\|_{L^2}
 \leq\|\mathbb{J}(w_N^{r-1}\psi)\|_{L^2}
      +d_2 \|w_N^r\psi\|_{L^2}.
      \label{eq:weighted-gradient-bound}
\end{align}
 Moreover, the weighted
interpolation inequality \cite[Lemma 1]{Ponce2011} gives
\begin{equation}\label{eq:weighted-interpolation}
 \|\mathbb{J}(w_N^{r-1}\psi)\|_{L^2}
 \leq d_3\|\mathbb{J}^{r}\psi\|_{L^2}^{1/r}
           \|w_N^r\psi\|_{L^2}^{(r-1)/r}
 =d_3\|\psi\|_{H^r}^{1/r}
           \|w_N^r\psi\|_{L^2}^{(r-1)/r},
\end{equation}
Combining \eqref{eq:linear-initial-bound},
\eqref{eq:weighted-gradient-bound} and
\eqref{eq:weighted-interpolation}, we arrive at
\begin{equation}\label{eq:linear-case-large}
 |I_N(t)|
 \leq d_4\|w_N^r\psi(t)\|_{L^2}
 \left(
   \|\psi(t)\|_{H^r}^{1/r}
      \|w_N^r\psi(t)\|_{L^2}^{(r-1)/r}
   +\|w_N^r\psi(t)\|_{L^2}
   +\|\psi(t)\|_{H^1}
 \right).
\end{equation}

It remains to estimate $II_N$. The Sobolev embedding
$H^1(\mathbb R)\hookrightarrow L^\infty(\mathbb R)$, the boundedness
of $\Pi$ on Sobolev spaces, and the fact that $H^2(\mathbb R)$ is an
algebra imply
\begin{align}
 \|\partial_x\Pi(\overline{u}\psi)\|_{L^\infty}
 \leq d_5 \|\partial_x\Pi(\overline{u}\psi)\|_{H^1}
 \leq d_6\|\partial_x(\overline{u}\psi)\|_{H^1}
 \leq d_7\|\overline{u}\psi\|_{H^2}
 \leq d_8\|u\|_{H^2}\|\psi\|_{H^2}.
      \label{eq:projected-product-bound}
\end{align}
Therefore, by H\"older's inequality and
\eqref{eq:projected-product-bound},
\begin{align}
 |II_N(t)|
 &\leq2\|w_N^r u(t)\|_{L^2}
           \|w_N^r\psi(t)\|_{L^2}
           \|\partial_x\Pi(\overline{u}\psi)(t)\|_{L^\infty}
           \notag\\
 &\leq d_9 \|u(t)\|_{H^2}\|\psi(t)\|_{H^2}
             \|w_N^r u(t)\|_{L^2}
             \|w_N^r\psi(t)\|_{L^2}.
             \label{eq:nonlinear-final-bound}
\end{align}
Combining \eqref{eq:linear-case-large} and \eqref{eq:nonlinear-final-bound}, we obtain 
\begin{equation}
\begin{split}
    \frac{d}{dt} \| w_N^r \psi(t) \|_{L^2}^2
    \le 
    &
   2 d_4\|w_N^r\psi(t)\|_{L^2}
 \left(
   \|\psi(t)\|_{H^r}^{1/r}
      \|w_N^r\psi(t)\|_{L^2}^{(r-1)/r}
   +\|w_N^r\psi(t)\|_{L^2}
   +\|\psi(t)\|_{H^1}
   \right)
   \\
   &
   +
    2 d_9 \|u(t)\|_{H^2}\|\psi(t)\|_{H^2}
             \|w_N^r u(t)\|_{L^2}
             \|w_N^r\psi(t)\|_{L^2}
\end{split}
\end{equation}
For $t\in [-\mathcal T, \mathcal T]$, by Theorem \ref{foc_invariant_HsLr}, suppose
\begin{equation}\label{eq:additional-bounds}
 S_{\mathcal{T}}:=\sup_{|t|\leq \mathcal T}
 \bigl(\|u(t)\|_{H^2}
       +\|u(t)\|_{L_r^2}
       +\|\psi(t)\|_{H^2}\bigr)<\infty.
\end{equation}

Using Young’s inequality $ab\le \frac{a^p}{p}+\frac{b^q}{q}$ with $a, b\ge 0$, $p,q>1$ and $\frac{1}{p}+\frac{1}{q}=1$, we infer that 
\begin{equation}
\begin{split}
     \frac{d}{dt} \| w_N^r \psi(t) \|_{L^2}^2
     &
    \le 
    \frac{2r-1}{2r} \| w_N^r \psi \|_{L^2}^2 + \frac{1}{2r}\times (2d_4)^{2r} \|\psi(t) \|_{H^r}^2
    +
     2 d_4\|w_N^r\psi(t)\|_{L^2}^2
     +
     \|w_N^r\psi(t)\|_{L^2}^2
     \\
     &
     \quad +
     d_4^2S_{\mathcal{T}}^2
     +
     \| w_N^r \psi(t) \|_{L^2}^2
     +
     d_9^2 S_{\mathcal{T}}^6
     \\
     &
     \le
     P_r \|
w_N^r \psi\|_{L^2}^2 + Q_{r, \mathcal{T}}    ,
\end{split}
\end{equation}
where 
\[
P_r= \frac{2r-1}{2r}+2d_4+2, 
\quad
Q_{r,\mathcal{T}}= \frac{1}{2r}\times (2d_4)^{2r} \sup_{ |t| \le \mathcal{T}} \| \psi(t)\|_{H^r}^2+d_9^2S_{\mathcal{T}}^6+ d_4^2S_{\mathcal{T}}^2<+\infty.
\]
Applying Gronwall's inequality in both time directions  and Fatou's  Lemma, the proof is complete.

\end{proof}

\begin{proof}[Proof of Theorem \ref{H^N_1}]
From Proposition \ref{pro-u-to-Lu} and  $u_0 \in H_r^m$, it follows that 
\begin{equation}
 \mathcal{L}_{u_0}(u_0), \mathcal{L}^2_{u_0}(u_0), \dots, \mathcal{L}^m_{u_0}(u_0)\in L_r^2.
\end{equation}
According to formula (4.25) in  \cite{SunLMP},
$u(t)=W(t)e^{-it\mathcal{L}^2_{u_0}}(u_0)\in L_r^2$. Recall the unitary equivalence formula $$\mathcal{L}_{u(t)}=W(t) \mathcal{L}_{u_0} W(t)^*: H^1\to L^2.$$ 

For $0\le k\le m$, set
$f_k:=\mathcal{L}_{u_0}^k u_0\in H^2$ and
\[
\psi_k(t):=W(t)e^{-it\mathcal{L}_{u_0}^2}f_k,
\]
Since $f_k\in H^2=\mathrm{Dom}(\mathcal{L}_{u_0}^2)$,
using
$\partial_tW(t)=\mathcal{B}_{u(t)}W(t)$, we obtain
\[
\begin{aligned}
\partial_t\psi_k(t)
&=\bigl(\mathcal{B}_{u(t)}
-i\mathcal{L}_{u(t)}^2\bigr)\psi_k(t)\\
&=i\partial_x^2\psi_k(t)
+2\sigma u(t)\partial_x
\Pi\bigl(\overline{u(t)}\psi_k(t)\bigr).
\end{aligned}
\]
Thus $\psi_k$ satisfies \eqref{eq:psi} with
$\psi_k(0)=\mathcal{L}_{u_0}^k u_0$.

For all $ 0 \le k\le m$, by Lemma \ref{W_exp_Lr2}, we have
\begin{equation}
\begin{split}
 \mathcal{L}^k_{u(t)}\big(u(t)\big)
 &
 =\mathcal{L}^k_{u(t)}(W(t)\exp(-it\mathcal{L}_{u_0}^2)(u_0))=W(t)\mathcal{L}^k_{u_0} (\exp(-it\mathcal{L}_{u_0}^2)(u_0))
 \\
 &
    =W(t)\exp(-it\mathcal{L}_{u_0}^2)\big( \mathcal{L}^k_{u_0}(u_0) \big)\in W(t)\exp(-it\mathcal{L}_{u_0}^2)\big( L_r^2\cap H^{p-k} \big) \in L_r^2.
\end{split}
\end{equation}
Finally, it follows from Proposition \ref{prop_Lu-to-u} that $u(t)\in H_r^m$. 

\end{proof}

\end{document}